\documentclass[11pt]{article}
\usepackage{amsmath,amssymb,amsthm}

\newcommand{\Tbar}{\ensuremath{\overline{\mathcal{T}}}}

\newcommand{\hess}{\operatorname{Hess}}
\newcommand{\grad}{\operatorname{grad}}
\newcommand{\tr}{\operatorname{tr}}

\newcommand{\dr}{\frac{\partial}{\partial r}}
\newcommand{\dthet}{\frac{\partial}{\partial \vartheta}}

\newcommand{\mathH}{\mathbb{H}}
\newcommand{\mathR}{\mathbb{R}}

\newcommand{\mathQ}{\mathbb{Q}}
\newcommand{\mathZ}{\mathbb{Z}}

\newcommand{\caB}{\mathcal{B}}
\newcommand{\caF}{\mathcal{F}}
\newcommand{\caG}{\mathcal{G}}

\newcommand{\caP}{\mathcal{P}}
\newcommand{\caS}{\mathcal{S}}
\newcommand{\caT}{\mathcal{T}}

\newtheorem{definition}{Definition}

\newtheorem{theorem}[definition]{Theorem}
\newtheorem{corollary}[definition]{Corollary}

\begin{document}

\title{The Weil-Petersson metric geometry\footnote{2000 Mathematics Subject Classification Primary: 32G15; Secondary: 20H10, 30F60.}}         
\author{Scott A. Wolpert}        % Enter your name between curly braces
\date{December 30, 2007}          % Enter your date or \today between curly braces
\maketitle
\begin{abstract}
A summary introduction of the Weil-Petersson metric space geometry is presented.  Teichm\"{u}ller space and its augmentation are described in terms of Fenchel-Nielsen coordinates.  Formulas for the gradients and Hessians of geodesic-length functions are presented. Applications are considered.  A description of the Weil-Petersson metric in Fenchel-Nielsen coordinates is presented.  The Alexandrov tangent cone at points of the augmentation is described.  A comparison dictionary is presented between the geometry of the space of flat tori and Teichm\"{u}ller space with the Weil-Petersson metric.  
\end{abstract}

\section{Introduction}       
Our goal is to present a summary introduction for the current understanding of a Weil-Petersson metric space.  With apologies to colleagues, our goal is to present an exposition following a development of concepts, rather than an exposition following the order of results discovered.  Selected readings and general attributions are provided at the end of each section.

There are overlapping themes for the current research on finite dimensional Weil-Petersson metrics.  Beginning with the work of Brock and in collaborations, the large scale coarse geometry is under extensive investigation \cite{BMi, Brkwp,Brkwpvs,BrMr,BrMs}.  Brock also initiated the consideration of the comparison to the geometry of quasi Fuchsian groups.    Beginning with the work of Yamada, the $CAT(0)$ geometry is also under examination \cite{DW2, MW, Wlcomp, Wlbhv, Yam2}.  Beginning with the work of Mirzakhani, the relationship to the Witten-Kontsevich conjecture and symplectic reduction are under continuing investigation \cite{Mirwitt, Mond, Saf, SafMu}.   The work of Mirzakhani combines explicit integrals and elements of Thurston's geometry to find the asymptotic count of lengths of simple closed geodesics on a hyperbolic surface \cite{Mirvol, Mirgrow}.   Following the work of Bridgeman and Taylor \cite{BrTr}, McMullen has shown that the metric can be reconstructed from dynamical quantities, such as measures on the unit circle and limit sets on the sphere \cite{McMthr}.  Beginning with the considerations of Weng \cite{Weng0,Weng}, the application to an arithmetic Riemann-Roch is being investigated \cite{GFM}.  McMullen \cite{McM}, the collaboration of Liu, Sun and Yau 
\cite{LSY1,LSY2}, as well as Yeung \cite{Yun, Yun1} have examined comparisons between the metric and the classical metrics for a domain and have also considered applications.  Huang continues a detailed examination of the curvature of the metric \cite{Zh,Zh2,Zh3,Zh4}.

\section{Basics of Teichm\"{u}ller theory}
Let $\caT$ be the Teichm\"{u}ller space for homotopy marked genus $g$, $n$ punctured Riemann surfaces $R$ of negative Euler characteristic.  A point of $\caT$ is the equivalence class of a pair $(R,f)$, with $f$ a homeomorphism from a reference topological surface $F$ to $R$.  By Uniformization a conformal structure determines a unique complete compatible hyperbolic metric $ds^2$ for $R$.  The Teichm\"{u}ller space is a complex manifold of dimension $3g-3+n$ with the cotangent space at $R$ represented by $Q(R)$, the space of holomorphic quadratic differentials on $R$ with at most simple poles at the punctures.  Weil introduced the Hermitian cometric. 
\begin{definition}
The Weil-Petersson cometric is $\langle\varphi,\psi\rangle\ = \ \int_R\varphi\overline\psi\,(ds^2)^{-1}.$
\end{definition}
In the 1940's Teichm\"{u}ller introduced the Finsler metric with conorm given as $\|\varphi\|_T=\int_R|\varphi|$.  
The Weil-Petersson (WP) dual metric is invariant under the action of the {\em mapping class group}, $MCG$, the group of homeomorphisms modulo the subgroup of homeomorphisms homotopic to the identity relative to punctures. The WP metric projects to the quotient $\mathcal M=\caT/MCG$, the moduli space of Riemann surfaces.  First properties are that the metric is K\"{a}hler, non complete with negative sectional curvature $\kappa$ with $\sup_\caT\kappa=0$ (except for $\dim\caT=1$ where 
$\sup_\caT\kappa<0$) and $\inf_\caT\kappa=-\infty$.  The metric continues to be the primary metric for understanding the K\"{a}hler geometry of Teichm\"{u}ller space \cite{McM, LSY1, LSY2, Wlpos, Yun1}.  The current exposition focuses on the metric space geometry.  {\em In practice and experience the WP geometry of $\caT$ corresponds to the hyperbolic geometry of surfaces.} 

A hyperbolic surface has a {\em thick-thin} decomposition with {\em thin} the region of injectivity radius below a threshold value.  The thin components of a hyperbolic surface are neighborhoods of cusps or are collars (fixed area tubular neighborhoods of short geodesics.)  Mumford first observed that the set of hyperbolic surfaces with lengths of closed geodesics bounded below by a constant $c>0$ forms a compact subset $\mathcal M_c$ of the moduli space $\mathcal M$. In general the totality of all thick regions of a given topological type forms a compact set of metric spaces in the Gromov-Hausdorff topology.  The Teichm\"{u}ller and WP metrics are comparable on $\mathcal M_c$.  The Teichm\"{u}ller and WP geometry of the ends 
of $\mathcal M-\mathcal M_c$ are examined in the references \cite{Minel,Msext,Wlcomp, Wlbhv}, as well as discussed below.  

We have for the reference topological surface $F$.
\begin{definition}
A $k$-simplex of the complex of curves $C(F)$ is a set of $k+1$ distinct free homotopy classes of non trivial, non peripheral, mutually disjoint simple closed curves of $F$. The pants graph $P(F)$ has vertices the maximal simplices 
of $C(F)$.  Vertices of $P(F)$ are connected by an edge provided the corresponding sets of free homotopy classes differ by replacing a single curve by a curve intersecting the original curve one or two times.  
\end{definition} 
The vertices of $C(F)$ are the free homotopy classes of non trivial, non peripheral simple closed curves.  The simplices of $C(F)$ are the convex sums of vertices.  The complex of curves $C(F)$ is a lattice, partially ordered by inclusion of simplices with maximal simplices, called {\em pants decompositions}, having dimension $3g-4+n$.   A pants decomposition decomposes a surface into a union of 
$2g-2+n$ three holed spheres.  The pants graph $P(F)$ becomes a metric space by specifying the edges to have unit-length.  Corresponding to a pants decomposition are global coordinates, Fenchel-Nielsen (FN) coordinates for $\caT$ given as gluing-parameters for constructing surfaces from right hyperbolic hexagons.  The construction begins with right hexagons which can be doubled across alternating edges to obtain a pair of pants, a genus zero hyperbolic surface with three geodesic boundaries with lengths free-parameters in $\mathR_{> 0}$.  Boundaries of pants of a common length can be abutted to construct a larger surface.  

A pants decomposition for $F$ provides a combinatorial scheme to abut boundaries of pants to obtain a hyperbolic surface of genus $g$ with $n$ punctures.  In abutting boundaries there is the free-parameter of the relative displacement of one boundary with respect to the other.  Overall for each abutting there are two free-parameters.  The first parameter is the common boundary geodesic-length $\ell$ valued in $\mathR_{>0}$.  The second parameter is the relative displacement $\tau$ valued in $\mathR$ measured in hyperbolic distance ($\tau$ is initially measured between appropriate footpoints and then analytically continued.)
\begin{theorem}
For a pants decomposition $\caP$ the FN coordinates $\Pi_{j\in\caP}(\ell_j,\tau_j):\rightarrow(\mathR_{>0}\times\mathR)^{3g-3+n}$ provide a real analytic equivalence for $\caT$.  The WP K\"{a}hler form is 
$\omega_{WP}\,=\,\frac12\sum_j d\ell_j\wedge d\tau_j$.
\end{theorem}
The FN coordinate expression for $\omega_{WP}$ is independent of the particular choice of pants decomposition.

The {\em augmented Teichm\"{u}ller space} $\Tbar$ is a partial compactification in the style of Bailey-Borel.  The space $\Tbar$ is important to the understanding of the Deligne-Mumford compactification of the moduli space $\mathcal M$ and for the WP geometry.  Frontier spaces are adjoined to $\caT$ corresponding to allowing geodesic-lengths $\ell_j$ to assume the value zero with the FN angles $\theta_j=2\pi\tau_j/\ell_j$ undefined (in polar coordinates the angle is undefined at the origin.)  The vanishing $\ell_j=0$ describes a degenerate hexagon with a side-length vanishing and the adjacent sides meeting at a common point on the circle at infinity for the hyperbolic plane.  The vanishing $\ell_j=0$ corresponds to a degenerate hyperbolic surface with a pair of cusps in place of a simple closed curve ($\gamma_j$ in the pants decomposition is now represented by the horocycles about the cusps.)  In general for a simplex $\sigma\subset C(F)$, the $\sigma$-null stratum is the space of structures $\mathcal S(\sigma)=\{R\ \mbox{degenerate}\mid\ell_{\alpha}(R)=0\mbox{ iff }\alpha\in\sigma\}$. The frontier spaces $\caF_{\caP}=\cup_{\sigma\subset\caP}\caS(\sigma)$ subordinate to a pants decomposition $\caP$ are adjoined to $\caT$ with a neighborhood basis for $\caT\cup\caF_{\caP}$ prescribed by the specification that
\[
((\ell_{\beta},\theta_{\beta}),\ell_{\alpha}):\caT\cup\caS(\sigma)\rightarrow\Pi_{\beta\notin\sigma}(\mathR_{>0}\times\mathR)\times\Pi_{\alpha\in\sigma}(\mathR_{\ge0})
\]
is a continuous map. For a simplex $\sigma\subset\caP,\,\caP'$, contained in two pants decompositions, the neighborhood systems are equivalent.  

A structural property is that the deformation spaces $\caS(\sigma)$ are products of lower dimensional Teichm\"{u}ller spaces and the limit of the tangential component of the WP metric of $\caT$ is simply the WP metric of the lower dimensional stratum $\caS(\sigma)$.  A sequence of marked hyperbolic surfaces $\{R_n\}$ converges in $\Tbar$ provided there is a simplex $\sigma$ contained in a pants decomposition $\mathcal P$ with $\ell_{\alpha}(R_n),\,\alpha\in\sigma$ limiting to zero (no convergence condition is placed on $\theta_{\alpha},\,\alpha\in\sigma$) and FN parameters $(\ell_{\beta},\theta_{\beta}),\,\beta\in\caP-\sigma$ converging. 
\begin{theorem}
The augmented Teichm\"{u}ller space $\Tbar=\caT\cup_{\sigma\in C(F)}\caS(\sigma)$ is a non locally compact stratified space.  The augmented Teichm\"{u}ller space is the WP completion of Teichm\"{u}ller space.
\end{theorem}
A point of $\Tbar-\caT$ represents a marked degenerate hyperbolic structure for which a simplex of $C(F)$ has each element represented by a pair of cusps.  The augmentation is also described as the Chabauty topology closure of the faithful cofinite representations of $\pi_1(R)$ into $PSL(2;\mathR)$ modulo conjugacies by $PSL(2;\mathR)$.  The quotient $\Tbar/MCG$ is topologically the Deligne-Mumford stable curve compactification of the moduli space of curves.  

We make the comparison between the upper half plane $\mathH$ as the space of marked flat tori and $\caT$ as the space of marked hyperbolic structures.  The comparison is explored in the following sections.  A point $z\in\mathH$ determines the lattice in $\mathbb C$ with basis vectors $\{1,z\}$.  A lattice change of basis is given by the action of the elliptic modular group $PSL(2;\mathZ)$.

\renewcommand{\arraystretch}{1.3}
\begin{tabular}[t]{|p{2.2 in}|p{2.2 in}|}
\hline
{\em flat structures} & {\em hyperbolic structures} \\
\hline
\hline
$\mathH$ the space of homotopy marked tori & $\caT$ the space of homotopy marked hyperbolic metrics \\
\hline
$PSL(2;\mathZ)$ & mapping class group $MCG$\\
\hline
$\mathH/PSL(2;\mathZ)$ moduli space of flat tori & $\mathcal M$ moduli space of Riemann surfaces \\
\hline
$\mathH\cup\mathbb Q$ with horoball topology & $\Tbar$ augmented Teichm\"{u}ller space \\
\hline
$\mathH\cup\mathbb Q/PSL(2;\mathZ)$ moduli space of stable elliptic curves & $\Tbar/MCG$ Deligne-Mumford moduli space of stable curves \\
\hline
horoballs $\{\Im A(z)\ge c\mid c < 1/2\}$, $A\in PSL(2;\mathZ)$ & Bers regions \\
\hline
Farey graph $\caG$ & complex of curves $C(F)$ and pants graph $P(F)$ \\
\hline
hyperbolic metric & WP metric (K\"{a}hler with negative curvature) \\
\hline
initial tangents to geodesics ending at $\mathbb Q$ are dense & initial tangents to geodesics ending at maximally degenerate structures are dense \\ 
\hline
for $A\in PSL(2;\mathZ)$, the function $-\log(\Im A(z))$ & for a closed geodesic $\alpha$, the root geodesic-length $\ell_{\alpha}^{1/2}$ \\
\hline
gradient $\mu_A=\grad\log(\Im A(z))$ with $\langle\mu_A,\mu_A\rangle =1$ & gradient $\lambda_{\alpha}=\grad\ell_{\alpha}^{1/2}$ with $\langle \lambda_{\alpha},\lambda_{\alpha}\rangle = 1/2\pi + O(\ell_{\alpha}^3)$ \\
\hline
for $A\in PSL(2;\mathZ)$, hyperbolic metric given as $(d\log(\Im A(z)))^2\,+\,(d\log(\Im A(z))\circ J)^2$ & 
for  a closed geodesic $\alpha$, WP metric given as $(d\ell_{\alpha}^{1/2})^2\,+\,(d\ell_{\alpha}^{1/2}\circ J)^2\,+\,O(\ell_{\alpha}^3)$ \\
\hline
$D_U(-\mu_A)=\langle J\mu_A,U\rangle J\mu_A$ & $D_U \lambda_{\alpha}=3\ell_{\alpha}^{1/2}\langle J\lambda_{\alpha},U\rangle\,J\lambda_{\alpha}+O(\ell_{\alpha}^{3/2})$  \\
\hline
$\hess(-\log(\Im A(z))\ge 0$ & $\hess \ell_{\alpha}>0,\ \hess \ell_{\alpha}^{1/2}>0$\\
\hline
\end{tabular}
\vspace{11 pt}

The Farey graph $\caG$ is realized in $\mathH\cup\mathQ$ by connecting vertices (rational numbers) $p/q$ and $r/s$ with a geodesic in $\mathH$ provided $|ps-qr|=1$.  Brock and Margalit examine geometric comparisons between the Farey graph and penultimate subsets of $C(F)$ \cite{BrMr} (the sets of decompositions containing a given $3g-5+n$ simplex of $C(F)$.)

The original reference for the WP metric is \cite{Ahsome}.  A reference for FN coordinates is \cite{Abbook} and for the symplectic form in FN coordinates is \cite{Wldtau}.   A counterpart approach to coordinates and the symplectic form for surfaces with cusps is extensively investigated in the works of Penner \cite{Pendec, Pencell} and also in \cite{Mond1}.  The complex of curves is introduced in \cite{Hrcomp} with a current introduction presented in \cite{HamH} and its metric space geometry is investigated in the foundational work of Masur and Minsky \cite{MaMiI, MaMiII}.  The original references for the augmented Teichm\"{u}ller space are \cite{Abdegn, Hrcomp} with the description in terms of the Chabauty closure of discrete faithful representations in \cite{HrCh}.  The original reference for analytic consideration of $\overline{\mathcal M}$ is \cite{Bersdeg} with \cite{Msext} presenting the first expansion for the metric. A brief survey of WP geometry current to the year 2002 is presented in the introduction of \cite{Wlcomp}.  Results from the work are also presented in the following sections.  Current understanding of WP curvature is presented in \cite{Zh, Zh2, Zh3}.

\section{The WP $CAT(0)$ geometry}
The augmented Teichm\"{u}ller space with the WP metric is a $CAT(0)$ metric space (a complete, simply connected, generalized non positively curved space.)  In particular $(\Tbar,d_{WP})$ is a length space, a metric space with unique distance-realizing paths (geodesics) between pairs of points.   Furthermore the WP metric has the Euclidean triangle comparison property: for a triangle in $\Tbar$ and a corresponding triangle in $\mathbb E^2$ with corresponding edge lengths, chords of the former (located by endpoints on sides) have lengths bounded by lengths of corresponding chords of the latter.  The $CAT(0)$ property is weaker than Gromov hyperbolicity, a generalized negative curvature property.  The strata $\caS(\sigma)\subset\Tbar,\,\sigma\in C(F)$ are intrinsic to the WP metric geometry and $\Tbar$ is an infinite polyhedron as follows.

\begin{theorem}
$\Tbar$ is a stratified metric space with each open strata characterized as the union of all geodesics containing a given point as an interior point.  The interior of a geodesic in $\Tbar$ is contained in a single stratum (geodesics do not refract at strata.)  $\Tbar$ itself is characterized as the closed convex hull of the maximally degenerate hyperbolic structures (the unions of thrice punctured spheres.)
\end{theorem}

The structure of strata provides that $MCG$ is the full group of orientation preserving isometries of $\caT$ as follows \cite{MW,Wlcomp}.   (The property does not follow the comparison between flat structures and hyperbolic structures, since the group of isometries of $\mathbb H$ is the Lie group $PSL(2;\mathbb R)$.)  A WP isometry extends to the completion $\Tbar$ and preserves the intrinsic strata structure, as well as the partial ordering of inclusion of simplices of $C(F)$.  Ivanov established that order preserving bijections of $C(F)$ are induced by elements of $MCG$ \cite{Ivaut}.  In particular for an isometry there is a corresponding element of $MCG$ and the two maps agree on the maximally degenerate structures in $\Tbar$ (a maximally degenerate structure is uniquely determined by its simplex.)  The two maps agree on $\Tbar$, the closed convex hull of the maximally degenerate structures. The isometry coincides with an element of $MCG$.      

There is a classification for the action of elements of $MCG$.  A mapping class $\iota$ acting on $\Tbar$ has fixed-points or positive translation length realized on a closed convex set $\mathcal A_{\iota}$, isometric to a metric space product $\mathR\times Y$.  In the latter case the isometry acts on $\mathR\times Y$ as the product of a translation and $id_Y$.  Following Thurston a mapping class is {\em irreducible} (pseudo Anosov) provided no power fixes the free homotopy class of a simple closed curve.  Infinite order non irreducible mapping classes are reducible, and are first analyzed in terms of mappings of proper subsurfaces.   Each irreducible mapping class has a unique invariant WP axis and non commuting irreducible mappings have divergent axes \cite{DW2,Wlcomp}.   

The pants graph $P(F)$ provides a quasi isometric model for $\caT$ and $\Tbar$ as follows.  Bers observed for constants $c'$ depending on the topological type $(g,n)$ that the sublevel sets 
$\caB_{\caP}=\{\ell_{\gamma}\le c'\mid\mbox{ for all }\gamma\in\caP \}$ for $\caP$ a pants decomposition, cover Teichm\"{u}ller space.  The bounded valence multivalued relation $ R\in\caT\leftrightarrow\{\caP\mid R\in\caB_{\caP}\}$ is the basis for (non unique) embeddings $h:\caT \rightarrow P(F)$ and $k:P(F)\rightarrow\caT$.  Brock established that the embeddings are quasi isometries (satisfy $d(x,y)/c'-c''\le d'(f(x),f(y))\le c'd(x,y)+c''$ for positive constants.)  On considering the edges of the pants graph $P(F)$ to have unit-length then Brock and in joint work with Margalit established the following \cite{Brkwp,BrMr}. 
\begin{theorem}
$\caT$ and $\Tbar$ are quasi isometric to $P(F)$.  For the topological types $(1,2)$ and $(0,4)$ the spaces $\caT$ and $\Tbar$ are quasi isometric to the Farey graph.
\end{theorem}   

Degenerate hyperbolic surfaces given as a union of surfaces of penultimate and ultimate types $(1,2)$, $(0,4)$ and $(0,3)$ are of special interest.  For such a simplex $\sigma\in C(F)$ with a total of $m$ subsurfaces of type $(1,2)$ or $(0,4)$ the corresponding strata $\caS(\sigma)\subset \Tbar$ is quasi isometric to an $m$-fold product of Farey graphs and contains a WP isometric image of $\mathR^m$ (a flat) as a product of geodesics from each of the cited factors.  More generally a quasi flat is a quasi isometric embedding of $\mathR^p$.  Quasi flats are important for understanding a geometry and are a tool for understanding quasi isometric rigidity in the setting of higher-rank symmetric spaces.  Behrstock and Minsky settled \cite{BMi} the open question on rank of $\caT$ showing that the maximal dimension of a quasi flat in $\caT$ is $\lfloor \frac{3g+n-2}{2}\rfloor$ (the maximal possible count of $(1,2)$ and $(0,4)$ subsurfaces.)  For $\dim\caT=3$ Brock and Masur have shown that any quasi flat is within a bounded distance of a strata quasi isometric to a product of Farey graphs \cite{BrMr}. 

In collaborations Behrstock-Minsky \cite{BMi} and Behrstock-Kleiner-Minsky-Mosher have been investigating the asymptotic cones $\mathcal{AC}$ of $MCG$ (the Gromov-Hausdorff limits of rescalings of the group word-metric.)  The main result of the first work is that the maximal dimension of a locally-compact subset of $\mathcal{AC}$ equals the maximal dimension of an Abelian subgroup of $MCG$.  In the second work rescaling limits in $\mathcal{AC}$ of flats and quasi flats of $MCG$, as well as the action on $\mathcal{AC}$ of quasi isometries are examined.  Hamenstaedt is also investigating the geometry of $MCG$ \cite{Ham}.  

Select readings for the section include the investigations of Brock \cite{Brkwp,Brkwpvs}, Daskalopoulos-Wentworth \cite{DW2} and the author \cite{Wlcomp}.   

\section{Geodesic-length functions}

Associated to each non trivial, non peripheral free homotopy class $\alpha$ on a marked hyperbolic surface is the length $\ell_{\alpha}$ of the unique geodesic in the free homotopy class.   Geodesic-lengths are explicit with $2 \cosh \ell_{\alpha}/2=\tr A$ for the free homotopy class $\alpha$ corresponding to the conjugacy class of $A$ in the deck group within $PSL(2;\mathR)$ and geodesic-lengths can be combined to provide coordinates for $\Tbar$.  Geodesic-length functions also have a direct relationship to WP geometry.  For a simple closed geodesic $\alpha$ on a hyperbolic surface, the surface can be {\em cut} along the geodesic to form two circle boundaries, which can then be identified by a relative rotation to form a new surface.  A flow on $\caT$ is defined by considering families of surfaces $\{R_t\}$ for which at time $t$ reference points from sides of the original geodesic are displaced by $t$ on $R_t$.  The infinitesimal generator, the FN vector field $t_{\alpha}$, and the WP gradient of geodesic-length satisfy the basic duality $2t_{\alpha}=J\grad\ell_{\alpha}$ for $J$ the almost complex structure of $\caT$. 

The relationship between hyperbolic geometry and WP geometry is displayed in the formulas for gradients.  The twist-length formula is
\[
\langle\grad\ell_{\alpha},J\grad\ell_{\beta}\rangle=4\,\omega_{WP}(t_{\alpha},t_{\beta})=-2\sum_{p\in\alpha\cap\beta}\cos\theta_p
\]
for the geodesics $\alpha,\,\beta$ and intersection angles $\theta_*$ on the hyperbolic surface.  The length-length formula for geodesics $\alpha,\beta$ with corresponding deck transformations $A,B,$ with corresponding axes $\tilde\alpha,\tilde\beta$ is 
\[
\langle\grad\ell_{\alpha},\grad\ell_{\beta}\rangle=\frac{2}{\pi}\bigl(\ell_{\alpha}\delta_{\alpha\beta}\ +{\ 
\sum}^{\prime}_{\langle A\rangle\backslash\Gamma/\langle B\rangle}(u\log\frac{u+1}{u-1}-2\bigr)\bigr)
\]
for the Kronecker delta $\delta_*$, where for $C\in\langle A\rangle\backslash\Gamma/\langle B\rangle$ then $u=u(\tilde\alpha,C(\tilde\beta))$ is the cosine of the intersection angle if the lifts $\tilde\alpha$ and $C(\tilde\beta)$ intersect and is otherwise $\cosh d(\tilde\alpha,C(\tilde\beta))$; for $\alpha=\beta$ the double-coset of the identity is omitted from the sum.  In the latter case the values $d(\tilde\alpha,C(\tilde\beta))$ are the lengths of the shortest representatives on the surface for the homotopy classes connecting $\alpha$ to $\beta$ with homotopy relative to $\alpha$ and relative to $\beta$.  The above formulas are the result of explicit integration of the Weil-Petersson product.

Select readings for the section are the investigations \cite{Mond, Rier, WlFN, Wlsymp}. 

\section{WP convexity and curvature}
The geodesic-length functions give rise to WP convex functions (functions strictly convex along geodesics.)  The length of a geodesic, the square root of length, the total length of a measured geodesic lamination, as well as for a surface with cusps the distance between unit-length horocycles are all WP strictly convex functions.  In particular the geodesic-length sublevel sets, as well as the strata of $\Tbar$ are convex sets.  Although a simple counterpart to the above gradient formulas is not yet available, there are bounds and expansions for the Riemannian Hessian (the intrinsic second derivative) of geodesic-length. 

Expansions for small geodesic-lengths for quantities on $\caT$ provide important information. The primary interest are quantities given as integrals on Riemann surfaces.  The WP metric, WP curvature, as well as the gradients and Hessians of geodesic-length are examples.  The approach for an expansion is based on understanding the integrand on the thick and thin regions of the surface.  Commonly the leading term of the expansion is the contribution of the collar zeroth rotational mode of the integrand with all other contributions higher order.  The expansion for the gradient and Hessian for small geodesic-lengths are examples.

\begin{theorem}
\label{grhe}
The variations of geodesic-length satisfy
\[
\langle\grad\ell_{\alpha}^{1/2},\grad\ell_{\beta}^{1/2}\rangle\ -\ \delta_{\alpha\beta}/2\pi\
\] 
is positive and bounded by $O(\ell_{\alpha}^{3/2}\ell_{\beta}^{3/2})$ and
\[
2\ell_{\alpha}\hess\ell_{\alpha}[U,U]\ -\ \dot\ell_{\alpha}^2[U]\ -\ 3\dot\ell_{\alpha}^2[JU]
\]
is positive and bounded by $O(\ell_{\alpha}^3\|U\|^2_{WP})$ for a tangent vector $U$ where for $c_0$ positive the remainder term constants are uniform for $\ell_{\alpha},\ell_{\beta}\le c_0$.
\end{theorem}
The Hessian is directly related to covariant differentiation by $\hess h[U,V]=\langle D_U\grad h,V\rangle$ for a smooth function $h$ and vector fields $U,V$.

\begin{corollary}
The WP connection $D$ is described for bounded geodesic-length, root gradient $\lambda_{\alpha}=\grad\ell_{\alpha}^{1/2}$, and a tangent vector $U$ by 
\[
D_U\lambda_{\alpha}\ =\ 3\ell_{\alpha}^{-1/2}\langle J\lambda_{\alpha},U\rangle J\lambda_{\alpha}\ +\ O(\ell_{\alpha}^{3/2}\|U\|_{WP}).
\]
\end{corollary}
A property of small geodesic-lengths and WP geodesics follows. For a geodesic $\gamma(t)$ with tangent field $\frac{d}{dt}$ the quantity $f(t)=\langle\lambda_{\alpha},\frac{d}{dt}\rangle^2+\langle J\lambda_{\alpha},\frac{d}{dt}\rangle^2$ has vanishing principal term for its first derivative.  The quantity $f(t)$ is Lipschitz along $\gamma(t)$ with constant $O(\ell_{\alpha}^{3/2})$.  
 
The estimates for WP sectional curvatures are also examples of small geodesic-length expansions.  The WP curvature of the span $\{\grad\ell_{\alpha},J\grad\ell_{\alpha}\}$ is $O(-\ell_{\alpha}^{-1})$.  Similarly for a pair of deformations approximately supported on different components of thick the corresponding curvature is $O(-\ell_{sys})$ for $\ell_{sys}$ the smallest geodesic-length.  For a pair of deformations approximately supported in the same component of thick the corresponding curvature is approximately the curvature for the limiting $2$-plane tangent to a stratum of $\Tbar$.

The section is based on the work \cite{Wlbhv} on behavior of geodesic-length and the works \cite{Zh, Zh2,Zh3} on the curvature of the metric.  

\section{Approaching degenerate hyperbolic structures}
A refined description of geodesic-length functions near degenerate hyperbolic structures provides further understanding of the WP metric.  We consider hyperbolic structures near a proper stratum $\caS(\sigma),\,\sigma\in C(F)$.  The closure of $\caS(\sigma)$ in $\Tbar$ is a union of strata $\cup_{\tau,\,\sigma\subseteq\tau}\caS(\tau)$. A closed convex subset $\overline{\caS(\sigma)}$ of a $CAT(0)$ space is the base of an orthogonal projection $\Pi_{\overline{\caS(\sigma)}}$ and we also consider the distance $d_{\overline{\caS(\sigma)}}$ to the stratum.  The projection is distance non increasing with fibers fibered by geodesics.  For the $k$-simplex $\sigma=\{\alpha_1,\dots,\alpha_{k+1}\}$ there is the overall bound on the distance of $R\in\Tbar$ to $\overline{\caS(\sigma)}$ 
\[
d_{\overline{\caS(\sigma)}}(R)\,\le\,(2\pi(\ell_{\alpha_1}(R)+\cdots+\ell_{\alpha_{k+1}}(R)))^{1/2}
\]
a consequence of the root geodesic-length convexity and the gradient pairing expansion in Theorem \ref{grhe}.  The bound displays the incompleteness of the metric.  The inequality compares to the formal equation for the hyperbolic plane $d(z,\infty)=-\log \Im z$ with the difference between the logarithm and the square root demonstrating the difference between the complete hyperbolic metric and the incomplete WP metric.  With the covering property of Bers regions it follows that for a constant depending only on topological type, each point of $\caT$ is within a fixed distance of a maximally degenerate structure.  

There is an approximation of a long WP geodesic segment $\stackrel{\frown}{pq}$ as follows.  At the ending point $q$ introduce a geodesic to a maximally degenerate structure $r$ and introduce a third geodesic from the beginning point $p$ to the maximally degenerate structure $r$.  The triangle $\Delta pqr$ has two long sides and a bounded length side. The comparison Euclidean triangle has a small angle at its beginning point.  A consequence of $CAT(0)$ is that the corresponding angle $\angle qpr$ in $\caT$ is likewise bounded, the desired approximation.  In particular for sufficiently long geodesics the corresponding angles are sufficiently small. An immediate consequence is that geodesics ending at maximally degenerate structures are dense in the space of geodesics.   

A description of geodesics ending at a point of $\Tbar-\caT$ is available.  A geodesic is {\em projecting} to $\overline{\caS(\sigma)}$ provided its projection is a point or equivalently it is length-minimizing to $\overline{\caS(\sigma)}$.  Projecting geodesics are almost described as integral curves of a constant sum of root gradient length functions.  In particular for a unit-speed projecting geodesic $\zeta$ to $\caS(\sigma)$ and the root gradients $\lambda_j=\grad\ell_{\alpha_j}$ there are constants $a_j$ such that the tangent field to the geodesic $\zeta$ satisfies
\[
\frac{d}{dt}\,=\, (2\pi)^{-1}\sum^{k+1}_{j=1}a_j\lambda_j\ + \ O(t^4)
\]
with $(2\pi)^{1/2}\|(a_j)\|_{Euclid}=1$ and the distance satisfies
\[
d_{\overline{\caS(\sigma)}}\, =\ (2\pi\sum^{k+1}_{j=1}\ell_{\alpha_j})^{1/2}\ +\ O(\sum^{k+1}_{j=1}\ell_{\alpha_j}^{5/2}).
\]
The distance formula prefigures the approximation that the WP metric in a neighborhood of $\overline{\caS(\sigma)}$ compares to the product of the (lower dimensional) WP metric on $\overline{\caS(\sigma)}$ and a universal metric for the normal bundle. The approximation is discussed in the next section. 

Selected readings for the section are \cite{Brkwpvs,Wlcomp,Wlbhv}. The basic reference for $CAT(0)$ geometry is \cite{BH}.
     
\section{Metrics and Fenchel-Nielsen coordinates}
 
Fenchel-Nielsen coordinates provide a straightforward description of hyperbolic surfaces and a parameterization for Teichm\"{u}ller space.  The WP metric does not provide the structure of a symmetric space (the full isometry group is discrete); the metric is not expected to have an elementary closed-form expression in FN coordinates.  Expansions and comparisons for the metric provide an alternative to an elementary expression.

The FN twist-length coordinates $(\ell_j,\tau_j)^{3g-3+n}_{j=1}$ for assembling hyperbolic pants provide global coordinates for $\mathcal T$ with the WP K\"{a}hler form 
\[
\omega_{WP}=\frac12\sum\,d\ell_j\wedge d\tau_j
\]
and on the Bers region $\{\ell_j<c'\}$ the metric comparisons 
\[
\langle\ ,\ \rangle \ \asymp\ \sum  (d\ell_j^{1/2})^2+(d\ell_j^{1/2}\circ J)^2\ \asymp\ 
\sum \hess\ell_j 
\]
for $J$ the almost complex structure with uniform comparability (given $c'>0$ there exist constants $c_1,\,c_2$ such that the metric is bounded above and below on the Bers region in terms of the constants and the given expressions.)   In a neighborhood of the maximally degenerate structure  $\{\ell_j=0\mid j=1\dots 3g-3\}$ the WP metric has the expansions 
\begin{align*}
\langle\ ,\ \rangle \ =&\ 2\pi\sum (d\ell_j^{1/2})^2+(d\ell_j^{1/2}\circ J)^2\ +\ O(\sum \ell_j^3\,\langle\ ,\ \rangle )\\ =&\ \frac{\pi}{6}\sum \frac{\hess\ell_j^2}{\ell_j} \ +\ O(\sum \ell_j^2\,\langle\ ,\ \rangle ).
\end{align*}
There are corresponding expansions for the neighborhood of a general stratum.  Metric incompleteness is immediate. 

There is a comparison between the Teichm\"{u}ller and WP geometry for the ends of the moduli space.  On a Bers region the Teichm\"{u}ller norm is comparable as follows
\[
\|\ \|_T^2 \ \asymp\ \sum \bigl((d\ell_j)^2+(d\ell_j\circ J)^2\bigr)\ell_j^{-2}.
\]
Each metric is comparable to a product of model metrics for the tangent planes $\bigcap_{k\ne j}\ker d\ell_k\cap\ker d\ell_k\circ J$.  The model for the Teichm\"{u}ller metric is the hyperbolic metric itself.  We now examine the model for the WP metric.  

There is a relation between the $1$-form $d\ell_j\circ J$ and the FN angle $\theta_j$.  The definition of the angle $\theta_j=2\pi\tau_j/\ell_j$ requires a pants decomposition.  An alternative definition of an angle is given by starting with the FN vector field $t_j=1/2\, J\grad \ell_j$, considering the WP dual, and defining the FN gauge $1$-form $\rho_j=2\pi(\ell_j^{3/2}\langle\lambda_j,\lambda_j\rangle)^{-1}\langle\ ,J \lambda\rangle$ for $\lambda_j=\grad \ell_j^{1/2}$.    The FN gauge is determined without the choice of a pants decomposition and satisfies the essential property $\rho_j(t_j)=d\theta_j(t_j)= 2\pi/\ell_j$.  Gauges and angles agree on the level sets of pants length $(\ell_j)_{j=1}^{3g-3}$.   

The WP metric for geodesic-lengths $\ell_j=2\pi^2r_j^2$ and FN gauges $\rho_j$ is comparable  as follows
\[
\langle\ ,\ \rangle \ \asymp\ \sum  4dr_j^2+r_j^6\varrho_j^2 
\]
with the expansion
\[
\langle\ ,\ \rangle \ = \pi^3 \sum  4dr_j^2+r_j^6\varrho_j^2 \ +\ O(\sum \ell_j^3 \,\langle\ ,\ \rangle ) 
\]
at the maximally degenerate structure.  There are corresponding expansions for the neighborhood of a general stratum.  The general expansion is in terms of a product of model metrics and the WP metric of the stratum.  

The model metric $4dr^2+r^6d\vartheta^2$  has K\"{a}hler form $2r^3drd\theta$,  Riemannian connection $D$ characterized by
\[
D_{\dr}\tfrac{\partial}{\partial r}=0,\quad D_{\dthet}\tfrac{\partial}{\partial r}=D_{\dr}\tfrac{\partial}{\partial\vartheta}=\tfrac3r\tfrac{\partial}{\partial\vartheta}\quad\mbox{and}\quad D_{\dthet}\tfrac{\partial}{\partial\vartheta}=\tfrac{-3}{4}r^5\tfrac{\partial}{\partial r}
\]
and Riemannian curvature $-3/2r^2$.  
The correspondence between the model and WP metrics is for $\frac{\partial}{\partial r_j}$ corresponding to $2^{3/2}\pi^2\lambda_j$ and $\frac{\partial}{\partial \vartheta_j}$ corresponding to the FN angle variation $(2\pi)^{-1}\ell_jt_j$.   The WP metric, K\"{a}hler form and connection have leading terms exactly corresponding to the expressions for the model metric. The WP curvature is comparable to the corresponding model curvature.  In effect the correspondence of metrics is a $C^1$ approximation with a bounded $C^2$ comparison.  

Readings for the section are \cite{Wlbhv, Wlext} with the comparability of the Teichm\"{u}ller metric presented in \cite{McM} and a refined large scale comparison presented 
in \cite{Minel}.

\section{WP Alexandrov tangent cone}
In a $CAT(0)$ metric space there is a well-defined angle between a pair of geodesics from a common initial point.  The angle enters in the definition of the tangent cone and in the first variation formula for distance.  At a point of the Teichm\"{u}ller space $\mathcal T$ the WP Alexandrov angle is given in terms of the Riemannian metric.  At a point of a stratum $\mathcal T(\sigma)\subset \Tbar -\mathcal T$, $\sigma\in C(F)$, the WP Alexandrov tangent cone $AC$ is isometric to a product of a Euclidean orthant and the tangent space $\mathbf T \mathcal T(\sigma)$ with the WP metric. The dimension of the Euclidean orthant is the count $|\sigma|$ of geodesic-lengths trivial on $\mathcal T(\sigma)$. 

A triple of points $(p,q,r)$ in $\Tbar$ has Euclidean comparison triangle with angle $\angle (p,q,r)$ valued in the interval $[0,\pi]$ determined by the Law of Cosines $2d(p,q)\,d(p,r)\cos\angle (p,q,r)=d(p,q)^2+d(p,r)^2-d(q,r)^2$.  The Alexandrov angle $(p,q,r)\rightarrow\angle (p,q,r)$ is upper semi continuous.   For constant speed geodesics $\gamma_0(t),\gamma_1(t)$ with common initial point  (from the $CAT(0)$ inequality) the comparison angle for  $(\gamma_0(0),\gamma_0(t),\gamma_1(t'))$ is a non decreasing function of $t$ and $t'$.  The Alexandrov angle is defined by the limit
\[
\cos \angle(\gamma_0,\gamma_1)=\lim_{t\rightarrow 0} \frac{d(p,\gamma_0(t))^2\,+\, d(p,\gamma_1(t))^2-d(\gamma_0(t),\gamma_1(t))^2}{2d(p,\gamma_0(t))\,d(p,\gamma_1(t))}.
\]
Geodesics at zero angle are said to define the {\em same direction}.  At zero angle provides an equivalence relation on the geodesics beginning at a point $p$ with the Alexandrov angle providing a metric on the space of directions.  The Alexandrov tangent cone $AC_p$ is the set of constant speed geodesics beginning at $p$ modulo the equivalence relation of same speed and at zero angle. 

A relative length basis for a point $p$ of $\mathcal T(\sigma)$ is a collection $\tau$ of vertices of $C(F)$  disjoint from the elements of $\sigma$ such that at $p$ the 
gradients $\{\grad \ell_{\beta}\}_{\beta\in\tau}$ span the tangent space $\mathbf T\mathcal T(\sigma)$.  A relative length basis can be given as the union of a partition and a dual 
partition for $R-\sigma$.  

We describe for the augmentation point $p$ an isometry between the Alexandrov tangent cone $AC_p$ and the product $\mathbb R_{\ge 0}^{|\sigma|}\times\mathbf T_p\mathcal T(\sigma)$ with the first factor the Euclidean orthant and the second factor the stratum tangent space with WP metric.   The mapping for a geodesic $\gamma(t)$ terminating at $p$ is given by associating for the lengths $\mathcal L(\gamma(t))=(\ell_{\alpha}^{1/2},\ell_{\beta}^{1/2})_{\alpha\in\sigma,\,\beta\in\tau}(\gamma(t))$  the initial one-sided derivatives  
\[
\Lambda:\gamma\rightarrow(2\pi)^{1/2}\frac{d\mathcal L(\gamma)}{dt}(0)
\]
(convexity provides for existence of the initial derivatives.)
By hypothesis the tuple $(\ell_{\beta}^{1/2})_{\beta\in\tau}$ provides local coordinates at $p$ for the stratum $\mathcal T(\sigma)$ and thus $(2\pi)^{1/2}\bigl(\frac{d\ell_{\beta}^{1/2}(\gamma)}{dt}(0)\bigr)_{\beta\in\tau}$ defines a vector in the tangent space $T_p\mathcal T(\sigma)$ with WP inner product.  The positive orthant $\mathbb R_{\ge 0}^{|\sigma|}\subset\mathbb R^{|\sigma|}$ is considered with the Euclidean inner product.  The Alexandrov tangent cone is given the structure of a cone in an inner product space through the formal relation $\langle\gamma_0,\gamma_1\rangle =\|\gamma_0'\|\|\gamma_1'\|\cos\angle (\gamma_0,\gamma_1)$.   
\begin{theorem}
The mapping $\Lambda$ from the WP Alexandrov tangent cone $AC_p$ to $\mathbb R_{\ge 0}^{|\sigma|}\times T_p\mathcal T(\sigma)$ is an isometry of cones with restrictions of inner products.   A WP terminating geodesic $\gamma$ with a root geodesic-length function initial derivative $\frac{d\ell_{\alpha}^{1/2}(\gamma)}{dt}(0)$ vanishing is contained entirely in the stratum $\{\ell_{\alpha}=0\}$.  Geodesics $\gamma_0$ and $\gamma_1$ at zero Alexandrov angle have comparison angles $\angle(p,\gamma_0(t),\gamma_1(t))$ bounded as $O(t)$.
\end{theorem}

A property for non positively curved Riemannian manifolds is that the exponential map is distance non decreasing.  An inverse exponential map $exp^{-1}_p:\Tbar\rightarrow AC_p$ is defined by associating to $q\in\Tbar$ the unique geodesic connecting $p$ to $q$ with speed $d(p,q)$.  The map is not an injection since geodesics at zero angle with common speed are mapped to a common element of $AC_p$.  The map is distance non increasing as follows.  From the $CAT(0)$ inequality and definition of the Alexandrov angle points on geodesics beginning at $p$ have distance satisfying
\begin{multline*}
d(\gamma_0(t),\gamma_1(t))^2 \ge d(p,\gamma_0(t))^2\,+d(p,\gamma_1(t))^2 \\ -2d(p,\gamma_0(t))\,d(p,\gamma_1(t))\cos\angle (\gamma_0,\gamma_1).
\end{multline*}
For equality for a single value the Flat Triangle Lemma provides that the geodesics are contained in a flat subspace of $\Tbar$.   The flat subspaces are classified.

A second application of the Alexandrov angle is the first variation formula for distance.  For the unit-speed geodesic $\gamma(t)$ the distance $d(\gamma(t),q)$ to a point not on the geodesic is convex with initial one-sided derivative satisfying
\[
\frac{d}{dt}d(\gamma(t),q)(0)=-\cos\angle (\gamma,\gamma_{pq})
\]
for $\gamma_{pq}$ the geodesic connecting $p$ to $q$.  {\em Non refraction} of geodesics on $\Tbar$ is a consequence: a WP length minimizing path at most changes strata at its endpoints. Consider a pair of unit-speed geodesics $\gamma_0(t),\, \gamma_1(t)$ with initial point $p$ such that the reverse path along $\gamma_0$ followed by $\gamma_1$ is length minimizing.  The Alexandrov angle between the tangents at $p$ is $\pi$.  The distance $d(\gamma_0(t),\gamma_1(t))$ is at least that of the path from $\gamma_0(t)$ to $p$ to $\gamma_1(t)$ and thus $\lim_{t\rightarrow 0} d(\gamma_0(t),\gamma_1(t))/2t=1$ and the angle is $\pi$.  Elements of $AC_p$ at angle $\pi$ necessarily lie in the subspace $\mathbf T_p\mathcal T(\sigma)$ and from the Theorem are segments of a single geodesic contained 
in $\mathcal T(\sigma)$. 

A third application is for length-minimizing paths connecting an initial and terminal point and intersecting a prescribed stratum.  Consider a pair of geodesics $\gamma_0$ and $\gamma_1$ each with initial point $p$ on $\overline{\mathcal T(\sigma)}$, $\gamma_0$ with endpoint $q$ and $\gamma_1$ with endpoint $r$.  Consider that the concatenation $\gamma_0+\gamma_1$ is a length-minimizing path connecting $q$ and $r$ to a point of $\overline{\mathcal T(\sigma)}$.  A geodesic $\kappa$ beginning at $p$ contained in $\overline{\mathcal T(\sigma)}$ provides a variation of the configuration.   The initial derivative of the distance $d(q,p)+d(r,p)$ along $\kappa$ is $-\cos\angle(\gamma_0,\kappa)-\cos\angle(\gamma_1,\kappa)$.  The geodesics beginning at $p$ contained in $\overline{\mathcal T(\sigma)}$ fill out the Alexandrov tangent cone $AC_p(\overline{\mathcal T(\sigma)})$.  It follows that the sum in $AC_p$ of the initial tangents of $\gamma_0$ and $\gamma_1$ has vanishing projection onto the subcone $AC_p(\overline{\mathcal T(\sigma)})$, the desired property.

A further application is for combinatorial harmonic maps.  Certain groups acting on Euclidean buildings and group extensions acting on Cayley graphs satisfying a Poincar\'{e} type inequality for links of points will have a global fixed point for an action on $\overline{\mathcal T}$.

The main reading for the section is \cite{Wlbhv}.  Additional readings are \cite{DW2,Wlcomp}.  The basic reference for Alexandrov angles is \cite{BH}.  The readings for combinatorial harmonic maps are \cite{IN,WgMT1,WgMT2}.   

%\nocite{WlO,PWWl}
%\bibliographystyle{alpha}
%\bibliography{scottbib}

\end{document}